\documentclass[12pt]{amsart}
\usepackage{amsmath,amsthm,latexsym,amscd,amsbsy,amssymb,amsfonts,fleqn,leqno,
euscript,pb-diagram,graphicx, texdraw}
\chardef\bslash=`\\ 

\makeatletter
\def\verbatim{\interlinepenalty\@M \@verbatim
  \leftskip\@totalleftmargin\advance\leftskip2pc
  \frenchspacing\@vobeyspaces \@xverbatim}
\makeatother
\makeatletter
  \def\dgt@k{\dg@DX=-3 \dg@DY=2 \dg@SIZE=3}
\makeatother

\makeatletter
  \def\dgt@kk{\dg@DX=3 \dg@DY=-1 \dg@SIZE=3}%
\makeatother

\theoremstyle{plain}
\newtheorem{thm}{Theorem}[section]
\newtheorem{cor}[thm]{Corollary}
\newtheorem{lem}[thm]{Lemma}
\newtheorem{pro}[thm]{Proposition}
\newtheorem{claim}[thm]{Claim}
\newcommand{\cl}{\mathrm{cl}}

\theoremstyle{definition}

\numberwithin{equation}{section}

\newcounter{rmnum}

\newenvironment{alphanum}{\begin{list}{{\rm (\alph{rmnum})}}
{\usecounter{rmnum}\def\makelabel##1{\hss\llap{##1}}
\setlength{\leftmargin}{0pt}\setlength{\itemindent}{37pt}
\setlength{\topsep}{5pt}\setlength{\parsep}{0pt}\setlength
{\itemsep}{0pt}}}{\end{list}}

\def\symbolnote#1#2{\let\thefootn=\thefootnote%
\renewcommand{\thefootnote}{\fnsymbol{footnote}}%
\footnotemark[#1]%
\footnotetext[#1]{#2}%
\let\thefootnote=\thefootn
}

\newfont{\bbb}{msbm10 scaled \magstep1}
\newfont{\bbc}{msbm8 scaled \magstep0}

\newcommand{\N}{\mbox{\bbb N}}

\begin{document}

\title{A non-separable Christensen's theorem and set tri-quotient maps}


\author{Stoyan  Nedev}
\address{Institute of Mathematics, Bulgarian Academy of Sciences,
Acad. G.Bonchev str., bl.8, Sofia 1113, Bulgaria}
\email{nedev@math.bas.bg}

\author{\boxed{\mbox{Jan Pelant}}}
\address{Mathematical Institute of Czech Academy of Sciences,
\v{Z}itn\'{a} 25, 11567 Praha 1, Czech Republic}

\author{Vesko Valov}
\address{Department of Computer Science and Mathematics, Nipissing University,
100 College Drive, P.O. Box 5002, North Bay, ON, P1B 8L7, Canada}
\email{veskov@nipissingu.ca}
\thanks{The last author was partially supported by NSERC
Grant 261914-03.}

\keywords{\v{C}ech completeness, set tri-quotient maps, sieve
completeness} \subjclass{Primary: 54C60; Secondary: 54E50.}

\begin{abstract}
For every space $X$ let $\mathcal K(X)$ be the set of all compact
subsets of $X$. Christensen \cite{c:74} proved that if $X, Y$ are
separable metrizable spaces and
$F\colon\mathcal{K}(X)\to\mathcal{K}(Y)$ is a monotone map such that
any $L\in\mathcal{K}(Y)$ is covered by $F(K)$ for some
$K\in\mathcal{K}(X)$, then $Y$ is complete provided $X$ is complete.
It is well known \cite{bgp} that this result is not true for
non-separable spaces. In this paper we discuss some additional
properties of $F$ which guarantee the validity of Christensen's
result for more general spaces.
\end{abstract}

\maketitle

\markboth{S.~Nedev, J.~Pelant, V.~Valov}{Set tri-quotient maps}

\section{Introduction}

All spaces in this paper are assumed to be completely regular.

The following characterization of Polish spaces established by J.P.
Christensen \cite{c:74} (see also \cite{r:73} for another proof) is
well known.

\begin{thm}\cite{c:74}
A separable metric space $Y$ is complete iff there exists a Polish
space $X$ and a map $F\colon\mathcal {K}(X)\to\mathcal {K}(Y)$ such
that:

\begin{itemize}
\item[(1)] $F$ is monotone $($i.e. if $K,L\in\mathcal{K}(X)$ with $K\subset L$, then $F(K)\subset F(L)$$)$;

\item[(2)] $F(\mathcal
{K}(X))$ is cofinal in $\mathcal{K}(Y)$ $($ i.e. for each $L\in
\mathcal{K}(Y)$ there is $K\in\mathcal{K}(X)$ with $L\subset
F(K)$$)$.
\end{itemize}
\end{thm}

According to Proposition 2.2(b) and Theorem 1.4 below, Theorem 1.1
remains valid if condition $(2)$ is replaced by the weaker one:

\begin{itemize}
\item[$(2)_c$] For any countable $L\in\mathcal{K}(Y)$ there exists
$K\in\mathcal{K}(X)$ with $L\subset F(K)$.
\end{itemize}

Theorem 1.1 is not valid for non-separable $X$. Indeed, let
$\mathbb{Q}$ be rational numbers and $X$ the discrete sum of all
compact subsets of  $\mathbb{Q}$. Then there exist a map $F:
{\mathcal K}(X)\rightarrow{\mathcal K}(\mathbb{Q})$ satisfying
conditions (1) and (2), see \cite{bgp}. Our first principal result
shows that Theorem 1.1 remains valid for arbitrary metrizable $X$
and $Y$ if $F$ satisfies an extra condition:

\begin{thm}
Let $X$ and $Y$ be metrizable spaces and  $F\colon\mathcal
{K}(X)\to\mathcal {K}(Y)$ be a map satisfying conditions $(1)$,
$(2)_c$ and condition $(3)_c$ below:
\begin{itemize}
\item[$(3)_c$] If $U\subset X$ and $V\subset Y$ are non-empty open sets
such that for each countable compact set $L\subset V$ there is a
compact $K\subset U$ with $L\subset F(K)$, then for any open cover
$\mathcal{W}$ of $U$ and any point $y\in V$ there exist a finite
subfamily $\mathcal{E}\subset\mathcal{W}$ and a neighborhood $V_y$
of $y$ such that each countable compact $K\subset V_y$ is covered by
$F(K)$ for some compact $K\subset\bigcup\mathcal{E}$.
\end{itemize}
Then $Y$ is completely metrizable and $densY\leq densX$ provided $X$
is completely metrizable.
\end{thm}

Any map $F\colon\mathcal {K}(X)\to\mathcal {K}(Y)$ satisfies $(3)_c$
if $X$ and $Y$ are metrizable with $X$ being separable (see
Proposition 2.2(b)). So, Theorem 1.2 is a generalization of
Christensen's result.

A non-metrizable analog of Theorem 1.1 was established in \cite{vd}
(see \cite{bc:96} for related results).

\begin{thm}\cite{vd}
Let $X$ be a Lindel\"{o}f \v{C}ech-complete space and
$F\colon\mathcal {K}(X)\to\mathcal {K}(Y)$ be a map satisfying
conditions $(1)$, $(2)$. If $Y$ is a $\mu$-complete $q$-space, then
$Y$ is also Lindel\"{o}f and \v{C}ech-complete.
\end{thm}

Recall that $X$ is said to be a {\em $\mu$-space} or {\em
$\mu$-complete} if every closed and bounded set in $X$ is compact.
Here, a set $A\subset X$ is {\em bounded} in $X$ if each continuous
real-valued function on $X$ is bounded on $A$. All paracompact, in
particular, Lindel\"{o}f spaces, are $\mu$-complete. The notion of a
$q$-space was introduced in \cite{michael:72}: $X$ is a $q$-space if
every $x\in X$ has a sequence $\{U_n\}$ of neighborhoods such that
if $x_n\in U_n$ for each $n\in\mathbb{N}$, then $\{x_n\}$ has a
cluster point in $X$. Obviously, every first countable, in
particular, every metric space is a $q$-space.

In order to obtain a general version  of Theorem 1.2 which implies
Theorem 1.3, we introduce a special type of set-valued maps called
{\em set tri-quotient maps} (see Section 2). Recall that
tri-quotient maps (single-valued) introduced by Michael
\cite{michael:77}  are extensively investigated, see \cite{jw1},
\cite{jw2}, \cite{m1}, \cite{m2}, \cite{o1}, \cite{p}, \cite{vu}.


Every map $F\colon\mathcal {K}(X)\to\mathcal {K}(Y)$ satisfying
conditions $(1)$, $(2)_c$ and $(3)_c$ is a monotone set tri-quotient
map (see Proposition 2.4). This allows us to derive Theorem 1.2 and
Theorem 1.3 from the following one which in turn follows from
Theorem 3.3 (recall that sieve-completeness, see \cite{ccn:74} and
\cite{michael:77}, is a more general property than \v
Cech-completeness and both they are equivalent in the class of
paracompact spaces).

\begin{thm}
Let $X$ be a sieve-complete space and $F: {\mathcal
K}(X)\rightarrow{\mathcal K}(Y)$ be a monotone set tri-quotient map.
If $Y$ is a $\mu$-space, then $Y$ is also seieve-complete and the
Lindel\"{o}f number $l(Y)$ of $Y$ is $\leq l(X)$.
\end{thm}

In the last section we apply Theorem 3.3 to show that sieve
completeness is preserved under linear continuous surjections
between function spaces, see Theorem 4.3. We also establish a
locally compact version of Theorem 1.2.

\section{Set tri-quotient maps}

The topology of a space $X$ is denoted by ${\mathcal T}(X)$.

Let ${\mathcal S}(X)\subset 2^{X}$. A map $F: {\mathcal S}(X)\rightarrow
2^{Y}$ is called {\em set tri-quotient\/} if there exists a map
$s: {\mathcal T}(X)\rightarrow {\mathcal T}(Y)$ such that:
\begin{itemize}
\item[(str1)] $s(U)\subset\bigcup\{F(K):K\in {\mathcal S}(X)\ \text{and}\ K\subset U\}$;
\item[(str2)] $s(X)=Y$;
\item[(str3)] $U\subset V$ implies $s(U)\subset s(V)$;
\item[(str4)] if $y\in s(U)$ and if ${\mathcal W}$ is a cover of
$\bigcup\{K\in F^{-1}(y): K\subset U\}$ by open subsets of $X$, then
$y\in s(\bigcup {\mathcal E})$ for some finite ${\mathcal E}\subset
{\mathcal W}$.
\end{itemize}

In the above definition $F^{-1}(y)$ stands for the family
$\{K\in{\mathcal S}(X): y\in F(K)\}$. Let us also observe that
conditions $(str1)$ and $(str4)$ imply that $F$ is surjective, i.e.
$Y=\bigcup\{F(K):K\in\mathcal{S}(X)\}$.

There is a similarity between set tri-quotient maps and Michael's
tri-quotient maps \cite{michael:77}. To clarify this similarity, let
us consider another class of maps introduced  in \cite{vd}.

A map $F: X\rightarrow 2^{Y}$ is said to be {\em generalized
tri-quotient\/} if one can assign to each open $U\subset X$ an open
$t(U)\subset Y$ such that:

\begin{itemize}
\item[(gtr1)] $t(U)\subset F(U)=\cup\{F(x): x\in U\}$;
\item[(gtr1)] $t(X)=Y$;
\item[(gtr1)] $U\subset V$ implies $t(U)\subset t(V)$;
\item[(gtr1)] if $y\in t(U)$ and if ${\mathcal W}$ is a cover
of $F^{-1}(y)\cap U$ by open subsets of $X$, then $y\in t(\bigcup
{\mathcal E})$ for some finite ${\mathcal E}\subset {\mathcal W}$.
\end{itemize}

We call the function $t: {\mathcal T}(X)\rightarrow {\mathcal T}(Y)$
an {\em assignment\/} for $F$. By (gtr1), every generalized
tri-quotient map is surjective, i.e. $Y=F(X)$. When $F: X\rightarrow
Y$ is single-valued and continuous, the above definition coincides
with the definition of a tri-quotient map \cite{michael:77}. It was
shown \cite[Proposition 2.1]{vd} that $F: X\rightarrow 2^{Y}$ is
generalized tri-quotient if and only if the projection $\pi _Y:
G(F)\rightarrow Y$ is tri-quotient, where $G(F)$ is the graph of
$F$.  This result, compared with \cite[Theorem 2.4]{o2}, shows that
generalized tri-quotient maps (as well as, set tri-quotient maps)
are different from the class of set-valued tri-quotient maps
introduced by Ostrovsky \cite{o2}.

Next lemma describes the connection between generalized
tri-quotient and set tri-quotient maps.

\begin{lem}
Let $F: X\rightarrow 2^Y$ be a generalized tri-quotient map. Then
$\Phi : 2^{X}\rightarrow 2^{Y}$, $\Phi (A)=\cl_Y(F(A))$, is monotone
set tri-quotient.
\end{lem}

\begin{proof}
It follows from the definition that $\Phi$ is monotone. Let $t:
{\mathcal T}(X)\rightarrow {\mathcal T}(Y)$ be an assignment for
$F$. We define $s(U)=t(U)$ for every open $U\subset X$. Obviously,
$s$ satisfies the first three conditions (str1)-(str3). Since
$F^{-1}(y)\cap U\subset\bigcup\{K\in\Phi ^{-1}(y): K\subset U\}$ for
all $y\in Y$ and $U\in{\mathcal T}(X)$, condition (str4) also holds.
\end{proof}

Similarly, every tri-quotient map $f: X\rightarrow Y$ generates a
monotone set tri-quotient map $F: {\mathcal
K}(X)\rightarrow{\mathcal K}(Y)$ defined by $F(K)=f(K)$,
$K\in{\mathcal K}(X)$.

Now, let us show that the map $F$ from Theorem 1.1 and Theorem 1.3
is monotone set tri-quotient.

\begin{pro}
Suppose $F: {\mathcal K}(X)\rightarrow{\mathcal K}(Y)$. Then we
have:
\begin{itemize} \item[(a)] $F$ is monotone set tri-quotient
provided $F$ satisfies conditions $(1)$ and $(2)$, $X$ is
Lindel\"{o}f and $Y$ a $\mu$-complete $q$-space;
\item[(b)] $F$ satisfies condition $(3)_c$ provided $X$ is separable metric and
$Y$ is first countable. Moreover, $F$ is monotone set tri-quotient
if $F$ satisfies conditions $(1)$ and $(2)_c$.
\end{itemize}
\end{pro}

\begin{proof}
To prove $(a)$, suppose $X$ is Lindel\"{o}f, $Y$ is a $\mu$-complete
$q$-space and $F$ satisfies conditions $(1)$ and $(2)$. We say that
a set $A\subset Y$ is $F$-covered by a set $B\subset X$ if for any
compact $L\subset A$ there exists a compact $K\subset B$ with
$L\subset F(K)$.

\begin{claim}
Let $U\subset X$ be functionally open and $V\subset Y$ open such
that $V$ is $F$-covered by $U$. If $\mathcal{W}$ is an open cover of
$U$ and $y\in V$, then there exists a neighborhood $V_y$ of $y$ and
a finite subfamily $\mathcal{E}\subset\mathcal{W}$ such that $V_y$
is $F$-covered by $\bigcup\mathcal{E}$.
\end{claim}

Since $U$ is functionally open, it is Lindel\"{o}f. So, we can
suppose that $\mathcal{W}=\{W_n:n\geq 1\}$ is countable. Let
$\{V_n\}$ be a sequence of neighborhoods of $y$ witnessing that $y$
is a $q$-point and such that $\cl(V_{n+1})\subset V_n\subset V$ for
all $n$. Assume the claim is false and for each $n$ choose a compact
set $L_n\subset V_n$ which is not covered by any $F(K)$,
$K\in\mathcal{K}\big(\bigcup_{i=1}^{i=n}W_i\big)$. Then the set

$$L=\big(\bigcup_{n=1}^{\infty}L_n\big)\bigcup\big(\bigcap
_{n=1}^{\infty}V_n\big)$$ is closed. It is bounded in $Y$ because
every infinite subset of $L$ has a cluster point. Hence $L$ is
compact (recall that $Y$ is a $\mu$-space). Since $L\subset V$ and
$V$ is $F$-covered by $U$, there is a compact set $K\subset U$ with
$L\subset F(K)$. Then $K\subset\bigcup_{i=1}^{i=m}W_i$ for some $m$.
Consequently, $L_m$ is covered by $F(K)$, which contradicts the
choice of $L_m$. The claim is proved.

Now, for every open $U\subset X$ let $s(U)$ be the set of all $y\in
Y$ having a neighborhood in $Y$ which is $F$-covered by a
functionally open subset $W$ of $X$ with $W\subset U$. Obviously,
$s(U)$ is open in $Y$ (possibly empty) and $s$ satisfies first three
conditions from the definition of a set tri-quotient map. To check
the last one, let $z\in s(U)$ and ${\mathcal W}$ be a cover of
$\bigcup\{K\in F^{-1}(z): K\subset U\}$ consisting of open in $X$
sets. Then there is a functionally open subset $W_0$ of $X$ with
$W_0\subset U$ and a neighborhood $V_0$ of $z$ such that $V_0$ is
$F$-covered by $W_0$. Since $F$ is monotone, $U=\bigcup\{K\in
F^{-1}(z): K\subset U\}$, so $\mathcal W$ is an open cover of $U$.
Taking a refinement of $\mathcal W$, if necessary, we can assume
that each element of $\mathcal W$ is functionally open in $X$. Then
$\mathcal{W}_0 =\{G\cap W_0: G\in{\mathcal W}\}$ is a functionally
open cover of $W_0$. According to Claim 2.3, there exist a
neighborhood $V_z$ of $z$ and finite
$\mathcal{E}_0\subset\mathcal{W}_0$ such that $V_z$ is $F$-covered
by $\bigcup\mathcal{E}_0$.

To finish the proof of $(a)$, let $\mathcal{E}=\{G\in\mathcal{W}:
G\cap W_0\in\mathcal{E}_0\}$. Because $V_z$ is $F$-covered by
$\bigcup\mathcal{E}_0$ which is functionally open in $X$ (as a
finite union of functionally open sets) and
$\bigcup\mathcal{E}_0\subset\bigcup\mathcal{E}$, we have that $z\in
s(\bigcup{\mathcal E})$. Therefore, $F$ is set tri-quotient and
monotone.

To prove $(b)$, assume $F$ does not satisfy $(3)_c$. Then there are
open sets $U\subset X$ and $V\subset Y$, an open cover $\mathcal{W}$
of $U$ and a point $y\in V$ such that every countable compact set
$L\subset V$ is covered by $F(K)$ for some compact set $K\subset U$,
but $y$ does not have a neighborhood which is contained in any
$\bigcup\{F(K): K\in\mathcal{K}(\bigcup\mathcal{E})\}$ with
$\mathcal{E}\subset\mathcal{W}$ being finite. Since $X$ is
separable, we can suppose $\mathcal{W}=\{W_n\}_{n\geq 1}$ is
countable. Next, choose neighborhoods $V_n\subset V$ of $Y$ and
countable compact sets $L_n\subset V_n$ such that $\{V_n\}_{n\geq
1}$ is a local base at $y$ and $L_n$ is not covered by any $F(K)$,
$K\subset\mathcal{K}\big(\bigcup_{i=1}^{i=n}W_i\big)$. Since
$L=\big(\bigcup_{n=1}^{\infty}L_n\big)\cup\{y\}$ is countable and
compact, there exists a compact set $K\subset U$ with $L\subset
F(K)$. As in the proof of Claim 2.3, this contradicts the choice of
the sets $L_n$. Hence, $F$ satisfies condition $(3)_c$.

It follows from Proposition 2.4 below that $F$ is monotone set
tri-quotient provided it satisfies conditions $(1)$ and $(2)_c$.
\end{proof}

\begin{pro}
Let $X$ and $Y$ be arbitrary spaces. Then any map $F: {\mathcal
K}(X)\rightarrow{\mathcal K}(Y)$ satisfying conditions $(1)$,
$(2)_c$ and $(3)_c$ is monotone set tri-quotient.
\end{pro}

\begin{proof}
Because $F$ satisfies $(1)$, it is monotone. For every open
$U\subset X$ we define $s(U)$ to be the set of all $y\in Y$ having a
neighborhood $V_y$ in $Y$ such that any countable compact $L\subset
V_y$ is covered by $F(K)$ for some compact set $K\subset U$.
Obviously, $s(U)$ is open in $Y$. Since $F$ satisfies conditions
$(1)$, $(2)_c$ and $(3)_c$, it is easily seen that $s$ satisfies
conditions $(str1)-(str4)$. So, $F$ is set tri-quotient.
\end{proof}

\section{Sieve-complete spaces}

\subsection{Proof of Theorem 1.4}
First, let us recall the definition of a sieve and a sieve-complete
space (see \cite{ccn:74} and \cite{michael:77}). A {\em sieve\/} on
a space $X$ is a sequence of open covers $\{U_{\alpha}: \alpha\in
A_n\}_{n\in\N}$ of $X$, together with maps $\pi _n:
A_{n+1}\rightarrow A_n$ such that $U_{\alpha}=\bigcup\{U_{\beta}:
\beta\in\pi ^{-1}_n(\alpha)\}$ for all $n$ and $\alpha\in A_n$. A
{\em $\pi$-chain\/} for such a sieve is a sequence $(\alpha _n)$
such that $\alpha _n\in A_n$ and $\pi (\alpha _{n+1})=\alpha _n$ for
all $n$. The sieve is {\em complete\/} if for every $\pi$-chain
$(\alpha_n)$, every filter base $\mathcal F$ on $X$ which meshes
with $\{U_{\alpha _n}: n\in \N\}$ (i.e. every $B\in\mathcal F$ meets
every $U_{\alpha _n}$) has a cluster point in $X$, or equivalently,
every filter base $\mathcal F$ on $X$ such that each $U_{\alpha _n}$
contains some $P\in{\mathcal F}$ clusters in $X$. A space $X$ with a
complete sieve is called {\em sieve-complete\/}. A sieve
$(\{U_{\alpha}: \alpha\in A_n\}, \pi _n)$ is said to be {\em
finitely additive\/} \cite{michael:77} if every cover $\{U_{\alpha}:
\alpha\in A_n\}$, as well as every collection of the form
$\{U_{\beta}: \beta\in\pi ^{-1}_n(\alpha)\}$ with $\alpha\in A_n$,
is closed under finite unions. When $\cl_X(U_{\beta})\subset
U_{\alpha}$ for all $\alpha\in A_n$ and $\beta\in\pi
^{-1}_n(\alpha)$, the sieve is called a {\em strong sieve\/}
\cite{ccn:74}.  Every sieve-complete space has a finitely additive
complete sieve \cite[Lemma 2.3]{michael:77}, as well as a strong
complete sieve \cite[Lemma 2.4]{michael:77}. Moreover, the proof of
\cite[Lemma 2.3]{michael:77} shows that the complete finitely
additive sieve which is obtained from a strong complete sieve is
also strong. Therefore, every sieve complete space has a strong
complete finitely additive sieve.

Let ${\mathcal S}(X)\subset 2^{X}$. We will use $\tau ^{+}_{V}$ to
denote the {\em upper Vietoris topology\/} on ${\mathcal S}(X)$
generated by all collections of the form $\hat{U}=\{H\in {\mathcal
S}(X): H\subset U\}$, where $U$ runs over the open subsets of $X$.

\begin{lem}
If $(\{U_{\alpha}: \alpha\in A_n\}, \pi _n)$ is finitely additive
and a strong complete sieve on $X$, then $(\{\hat{U}_{\alpha}:
\alpha
\in A_n\}, \pi _n)$ is a complete sieve on $({\mathcal K}(X), \tau ^{+}_{V})$.
\end{lem}

\begin{proof}
Because $\gamma=(\{U_{\alpha}: \alpha\in A_n\}, \pi _n)$ is a
finitely additive sieve on $X$, $\hat{\gamma} =(\{\hat{U}_{\alpha}:
\alpha\in A_n\}, \pi _n)$ is a sieve on $({\mathcal K}(X), \tau
^{+}_{V})$. Let us show that $\hat{\gamma}$ is complete. Suppose
$(\alpha _n)$ is a $\pi$-chain and $\mathcal F$ a filter base on
${\mathcal K}(X)$ which meshes with $\{\hat{U}_{\alpha _n}\}$. By
\cite[Lemma 2.5]{michael:77}, $K=\bigcap U_{\alpha _n}$ is a
nonempty compact subset of $X$ such that every open $W\supset K$
contains some $U_{\alpha _n}$. Then every neighborhood $\hat{W}$ of
$K$ in $({\mathcal K}(X), \tau ^{+}_{V})$ contains some
$\hat{U}_{\alpha_n}$, hence $\hat{W}$ meets every $H\in{\mathcal
F}$. Therefore $K$ belongs to the closure (in $({\mathcal K}(X),
\tau ^{+}_{V})$) of each $H\in{\mathcal F}$, i.e. $K$ is a cluster
point of ${\mathcal F}$ in $({\mathcal K}(X), \tau^{+}_{V})$.
\end{proof}

The following analogue of $q$-spaces was introduced in \cite{vv} :
call $X$ a {\it $wq$-space} if every $x\in X$ has a sequence
$\left\{U_n\right\}$ of neighborhoods such that if $x_n\in U_n$ for
each $n$, then $\left\{x_n\right\}$ is bounded in $X$. The
$wq$-space property is weaker than $q$-space property and they are
equivalent for $\mu$-spaces.\medskip

We say that a set-valued map $F: X\rightarrow 2^Y$ is a {\it
$wq$-map} if every $x\in X$ has a sequence $\{U_n\}$ of
neighborhoods such that if $x_n\in U_n$ for each $n$, then
$\bigcup\{F(x_n): n\in\N\}$ has a compact closure in $Y$. A version
of next lemma was established first in \cite[Lemma 2.3]{vd}. In the
present form it appears in \cite[Proposition 3.14]{vv}, and later on
in \cite[Theorem 2.2]{ch}.

\begin{lem}
\cite{vv} Let $F: X\rightarrow 2^Y$ be a $wq$-map with $Y$ being a
$\mu$-space. Then there exists an usco map $\Phi :X\rightarrow Y$
such that $F(x)\subset\Phi (x)$ for every $x\in X$.
\end{lem}

Next theorem provides the proof of Theorem 1.4.
\begin{thm}
Let $X$ be a  sieve-complete space and $Y$ a $\mu$-space. If there
exists a monotone set tri-quotient map $F: {\mathcal
K}(X)\rightarrow 2^Y$ such that each $F(K)$, $K\in{\mathcal K}(X)$,
has a compact closure in $Y$, then $Y$ is sieve-complete and
$l(Y)\leq l(X)$.
\end{thm}

\begin{proof}
As we already mentioned, there exists a strong complete sieve
$\gamma =(\{U_{\alpha}: {\alpha}\in A_n\}, \pi _n)$ on $X$ which is
finitely additive. Then, according to Lemma 3.1, $\hat{\gamma}$ is a
complete sieve on $({\mathcal K}(X), \tau ^{+}_{V})$.

First, let us show that $F$, considered as a set-valued map from
$({\mathcal K}(X), \tau ^{+}_{V})$ into $Y$, is a $wq$-map. Since
$\gamma$ is a finitely additive and strong sieve on $X$, for every
$K\in{\mathcal K}(X)$ there is a chain $(\alpha _n)$ such that
$K\subset U_{\alpha _n}$ for all $n$. This yields (see \cite[Lemma
2.5]{michael:77}) that $C=\bigcap U_{\alpha _n}$ is compact and
$\{U_{\alpha _n}\}$ is a base for $C$. We assign to $K$ the sequence
$\{\hat{U}_{\alpha _n}\}$. If $K_n\in\hat{U}_{\alpha _n}$ for all
$n$, then $H=(\bigcup K_n)\cup C$ is a compact subset of $X$ and,
since $F$ is monotone, $\bigcup F(K_n)\subset F(H)$. So, $\bigcup
F(K_n)$ has a compact closure in $Y$. Therefore $F$ is a $wq$-map
and, by Lemma 3.2, there exists an usco map $\Phi :({\mathcal K}(X),
\tau ^{+}_{V})\rightarrow Y$ with $F(K)\subset\Phi (K)$ for every
$K\in{\mathcal K}(X)$. Let us observe that $\Phi$ is onto, i.e.
$Y=\bigcup\{\Phi(K):K\in\mathcal{K}(X)\}$.  Since the Lindel\"{o}f
number of $({\mathcal K}(X), \tau ^{+}_{V})$ is $\leq l(X)$, the
last equality yields $l(Y)\leq l(X)$.

Because $F$ is set tri-quotient, there is a map $s: {\mathcal
T}(X)\rightarrow {\mathcal T}(Y)$ satisfying conditions
(str1)-(str4). Let $W_{\alpha}=s(U_{\alpha})$ for every $n$ and
$\alpha\in A_n$. We are going to show that $\lambda =(\{W_{\alpha}:
\alpha\in A_n\}, \pi _n)$ is a complete sieve on $Y$. Since all
$\gamma _n=\{U_{\alpha}: \alpha\in A_n\}$ are open covers of $X$, it
follows from conditions (str2) and (str4) that each $y\in Y$ is
contained in $s(\bigcup\omega _n)$ for some finite $\omega
_n\subset\gamma _n$. But each $\gamma _n$ is finitely additive, so
all systems $\{W_{\alpha}: \alpha\in A_n\}$, $n\geq 1$, are covers
of $Y$. Similarly, we can show that
$W_{\alpha}\subset\bigcup\{W_{\beta}: \beta\in\pi ^{-1}_n(\alpha)\}$
for every $n$ and $\alpha\in A_n$. The inclusions
$\bigcup\{W_{\beta}: \beta\in\pi ^{-1}_n(\alpha)\}\subset
W_{\alpha}$ follow from (str3) and $U_{\alpha}=\bigcup\{U_{\beta}:
\beta\in\pi ^{-1}_n(\alpha)\}$. Therefore, $\lambda$ is a sieve on
$Y$. To show that that $\lambda$ is a complete sieve, suppose
$(\alpha _n)$ is a $\pi$-chain and $\mathcal F$ is a filter base on
$Y$ which meshes with $\displaystyle \{W_{\alpha _n}: n\in\N\}$.
Then $\Phi ^{-1}({\mathcal F})=\{\Phi ^{-1}(P): P\in{\mathcal F}\}$
is a filter base on $({\mathcal K}(X), \tau ^{+}_{V})$.

\begin{claim}
$\Phi ^{-1}(\mathcal F)$ meshes with $\displaystyle\{\hat{U}_{\alpha
_n}: n\in\N\}$.
\end{claim}

If $\displaystyle y\in P\cap W_{\alpha _n}$ for some $P\in\mathcal
F$ and $n\in\N$, then, by (str1), there is $K\in{\mathcal K}(X)$
with $\displaystyle K\subset U_{\alpha _n}$ and $y\in
F(K)\subset\Phi (K)$. Therefore, $\displaystyle K\in\Phi
^{-1}(P)\cap\hat{U}_{\alpha _n}$ which completes the proof of the
claim.

Since $\hat{\gamma}$ is a complete sieve, $\Phi ^{-1}(\mathcal F)$ has a cluster
point, say $K_0$, in $({\mathcal K}(X), \tau ^{+}_{V})$.

\begin{claim}
$\Phi (K_0)\cap \cl_Y(P)\neq\varnothing$ for each $P\in{\mathcal
F}$.
\end{claim}

Suppose $\Phi (K_0)\cap \cl_Y(P)=\varnothing$ for some $P\in\mathcal
F$. Let $V\subset Y$ be open, disjoint with $P$ and containing $\Phi
(K_0)$. Because $\Phi$ is usc, there is a neighborhood $\hat{U}$ of
$K_0$ in $({\mathcal K}(X), \tau ^{+}_{V})$ such that $\Phi
(K)\subset V$ for every $K\in\hat{U}$. Since $\hat{U}$ meets $\Phi
^{-1}(P)$, $\Phi (K)\subset V$ for some $K\in\Phi ^{-1}(P)$ which is
a contradiction.

By Claim 3.5, ${\mathcal F}_0=\{\Phi (K_0)\cap \cl_Y(P):
P\in\mathcal F\}$ is a filter base on $\Phi (K_0)$. Because $\Phi
(K_0)$ is compact, ${\mathcal F}_0$ has a cluster point. So,
$\mathcal F$ has a cluster point in $Y$ and $\lambda$ is a complete
sieve on $Y$.
\end{proof}

Let us observe that the restriction in Theorem 3.3 $Y$ to be
$\mu$-complete and $F$ to be monotone were used only to apply Lemma
3.2 in order to find an usco map $\Phi : ({\mathcal K}(X),\tau
^{+}_V)\rightarrow{\mathcal K}(Y)$ with $F(K)\subset\Phi (K)$,
$K\in{\mathcal K}(X)$. Therefore, the following statement holds:

\begin{cor}
Let $F: ({\mathcal K}(X),\tau ^{+}_V)\rightarrow {\mathcal K}(Y)$ be
usc and set tri-quotient with $X$ a sieve-complete space. Then $Y$
is also sieve-complete.
\end{cor}

\begin{cor}
For a $\mu$-space $Y$ the following are equivalent:

\begin{alphanum}
\item $Y$ is sieve-complete.

\item There exists a paracompact \v Cech-complete space $X$ and
an open (not necessary continuous) surjection $f: X\rightarrow Y$
such that $f(K)$ has a compact closure in $Y$ for every $K\in {\mathcal
K}(X)$.

\item There exists a paracompact \v Cech-complete space $X$ and
a monotone set tri-quotient map $F: {\mathcal K}(X)\rightarrow {\mathcal K}(Y)$.
\end{alphanum}
\end{cor}

\begin{proof}
(a)$\ \Rightarrow\ $(b). This implication follows from \cite[Theorem
3.7]{ccn:74} stating that every sieve-complete space is an open and
continuous image of a paracompact \v Cech-complete space.

(b)$\ \Rightarrow\ $(c). If $f$ satisfies (b), we simply define $F:
{\mathcal K}(X)\rightarrow {\mathcal K}(Y)$ by $F(K)=\cl_Yf(K)$. Since $f$ is open,
$F$ is set tri-quotient.

(c)$\ \Rightarrow\ $(a). This implication follows from Theorem 3.3.
\end{proof}

\subsection{Proof of Theorem 1.2}
According to Proposition 2.4, Theorem 3.3 and the fact that sieve
and \v Cech-completeness are equivalent in the realm of paracompact
spaces, it follows that $Y$ is complete. Moreover, Theorem 3.3 also
implies that $dens Y\leq dens X$.

\section{Remarks and some applications}
Let us consider the following analogs of condition $(3)_c$ in
Theorem 1.2:
\begin{itemize}
\item[$(3)$] If $U\subset X$ and $V\subset Y$ are non-empty open sets
such that for each compact $L\subset V$ there is a compact $K\subset
U$ with $L\subset F(K)$, then for any open cover $\mathcal{W}$ of
$U$ and any point $y\in V$ there exists a finite subfamily
$\mathcal{E}\subset\mathcal{W}$ and a neighborhood $V_y$ of $y$ such
that for each compact $L\subset V_y$ there is a compact
$K\subset\bigcup\mathcal{E}$ with $L\subset F(K)$.
\item[$(3')$] For each open cover $\mathcal{W}$ of $X$
and for each point $y\in Y$ there exists a finite subfamily
$\mathcal{E}\subset\mathcal{W}$ and a neighborhood $V_y$ such that
every compact $L\subset V_y$ is covered by $F(K)$ for some compact
$K\subset\bigcup\mathcal{E}$.
\end{itemize}
Obviously, conditions $(3)$ and $(3)_c$ are not comparable, while
conditions $(2)$ and $(3)$ imply $(3')$. As in Lemma 2.2(b), one can
show that any map $F\colon\mathcal{K}(X)\to\mathcal{K}(Y)$ satisfies
condition $(3)$ if $X$ is second countable and $Y$ first countable.
Moreover, we have the following lemma whose proof is similar to that
one of Proposition 2.4.

\begin{lem}
If $X$ and $Y$ are arbitrary spaces and
$F\colon\mathcal{K}(X)\to\mathcal{K}(Y)$ satisfies conditions $(1)$,
$(2)$ and $(3)$, then $F$ is monotone set tri-quotient.
\end{lem}

We do not know whether Theorem 1.2 is valid when $F$ satisfies
conditions $(1)$, $(2)$ and $(3')$. It seems now that the related
claim in \cite[Theorem 5.2]{bgp} was overoptimistic.

It is interesting that a locally compact version of Theorem 1.2 is
true if $F$ satisfies conditions $(1)$ and $(3')$.

\begin{pro}
Let $X$ be a locally compact space and $F\colon\mathcal
{K}(X)\to\mathcal {K}(Y)$ satisfy conditions $(1)$ and $(3')$. Then
$Y$ is also locally compact.
\end{pro}

\begin{proof}
Let $\mathcal{U}=\{U_\alpha:\alpha\in A\}$ be an open cover of $X$
such that each $U_\alpha$ has a compact closure in $X$. Since $F$
satisfies condition $(3')$, for every $y\in Y$ there exists a
neighborhood $V_y$ and a finite $\mathcal{E}_y\subset\mathcal{U}$
such that every compact set $L\subset V_y$ is covered by $F(K)$ for
some compact $K\subset\bigcup\mathcal{E}_y$. So,
$V_y\subset\bigcup\{F(K):K\in\mathcal{K}(U_y)\}$, where
$U_y=\bigcup\{U:U\in\mathcal{E}_y\}$. Because the closure
$\overline{U}_y$ is compact and $F$ is monotone,
$\bigcup\{F(K):K\in\mathcal{K}(U_y)\}\subset
F\big(\overline{U}_y\big)$. Hence, each $V_y$ has a compact closure
in $Y$.
\end{proof}

As we already observed, if $X$ is second countable and $Y$ first
countable, then condition $(2)$ implies condition $(3')$. In this
case, Proposition 4.2 is valid whenever $F$ satisfies conditions
$(1)$ and $(2)$. The example provided in the introduction shows that
conditions $(1)$ and $(2)$ are not enough for the validity of
Proposition 4.2 if $X$ is not separable.

We are going now to apply Theorem 3.3 for obtaining alternative
proofs and improvements of some results from \cite{bgp} and
\cite{vv} concerning preservation of \v{C}ech completeness under
linear surjections between function spaces. Everywhere below
$C(X,E)$ denotes the set of all continuous maps from $X$ into $E$
(we write $C_p(X)$ when consider real-valued functions). The set
$C(X,E)$ endowed with the compact-open or the pointwise convergence
topology is denoted by $C_k(X,E)$ or $C_p(X,E)$, respectively. If
$u\colon C_k(X,E)\to C_p(Y,F)$ is a linear map, where $E$ and $F$
are normed spaces, then for every $y\in Y$ there exists a continuous
linear map $\mu_y\colon C_k(X,E)\to F$ defined by
$\mu_y(f)=u(f)(y)$, $f\in C_k(X,E)$. Following Arhangel'skii
\cite{ar}, we define the support $supp\big(\mu_y\big)$ of $\mu_y$ to
be the set of all $x\in X$ such that for every neighborhood $U$ of
$x$ in $X$ there is $f\in C(X,E)$ with $f(X\backslash U)=0$ and
$\mu_y(f)\neq 0$, see \cite{vv}. So, we can consider the set-valued
map $\varphi\colon Y\to 2^X$, $\varphi(y)=supp\big(\mu_y\big)$. This
map has the following properties (see \cite{bg}, \cite{vv}):
\begin{itemize}
\item[(a)] $\varphi$ is lower semi-continuous;
\item[(b)] if $L$ is a bounded set in $Y$, then $\varphi(L)$
is bounded in $X$;
\item[(c)] if $K$ is a bounded set in $X$, then the set $\varphi^{*}(K)=\{y\in
Y:\varphi(y)\subset K\}$ is bounded in $Y$;
\item[(d)] if $u$ is surjective, then $\varphi(y)\neq\varnothing$ for
all $y\in Y$.
\end{itemize}

It is shown in \cite[Theorem 3.3]{bgp} that if $u\colon C_p(X)\to
C_p(Y)$ is a continuous linear surjection with $X$ and $Y$
metrizable, then $Y$ is \v{C}ech-complete provided so is $X$. This
result was generalized in \cite[Corollary 3.15]{vv} to the case of
non-metrizable $X$ and $Y$ and function spaces of maps into normed
spaces (see the hypotheses of Theorem 4.3 below). Under the same
hypotheses, we can establish a sieve completeness version of this
result. Of course, if $X$ and $Y$ are paracompact spaces, then
Theorem 4.3 and \cite[Corollary 3.15]{vv} are equivalent. In such a
situation, Theorem 4.3 provides an alternative proof of
\cite[Corollary 3.15]{vv}.

\begin{thm}
Let $u\colon C_k(X,E)\to C_p(Y,F)$ be a continuous linear
surjection, where both $X$ and $Y$ are $\mu$-spaces and $Y$ a
$wq$-space. If $X$ is sieve-complete, then so is $Y$.
\end{thm}

\begin{proof}
Since $X$ is a $\mu$-space, $Y$ is a $wq$-space and $\varphi$
satisfies condition $(b)$, $\varphi$ is a $wq$-map. So, by Lemma
3.2, there exists an usco map $\phi\colon Y\to 2^X$ such that
$\varphi(y)\subset\phi(y)$ for every $y\in Y$. Now, define the  map
$F\colon\mathcal{K}(X)\to 2^Y$ by $F(K)=\phi^{*}(K)$. Let us note
that $F(K)$ may not be a compact subset of $Y$, but it has a compact
closure in $Y$. Indeed, $F(K)\subset\varphi^*(K)$ and the
$\mu$-completeness of $Y$ implies that the set $\varphi^*(K)$ is
compact as a closed and bounded subset of $Y$ (it is closed because
$\varphi$ is lower semi-continuous, and it is bounded because of
$(c)$). For every open $U\subset$ let $s(U)=\phi^*(U)$. Since $\phi$
is upper semi-continuous, every $s(U)$ is open in $Y$. We are going
to show that $s$ satisfies conditions $(str1)-(str4)$. Because
$\varphi(y)\neq\varnothing$ for all $y\in Y$, the sets $\phi(y)$,
$y\in Y$, are non-empty and compact. This yields that $s$ satisfies
conditions $(str1)$ and $(str2)$. Obviously, condition $(str3)$ also
holds. Finally, if $y\in s(U)$ and $\mathcal{W}$ is an open cover of
$U$, then $\phi(y)\subset U$ and choose a finite family
$\mathcal{E}\subset\mathcal{W}$ covering $\phi(y)$. So, $y\in
s\big(\bigcup\mathcal{E}\big)$. Therefore, $F$ is set tri-quotient
and we can apply Theorem 3.3 to conclude that $Y$ is sieve-complete.
\end{proof}

\ifx\undefined\bysame
\newcommand{\bysame}{\leavevmode\hbox to3em{\hrulefill}\,}
\fi

\end{document}